\documentclass{gtart}
\input gtmonout
\volumenumber{1}
\volumeyear{1998}
\volumename{The Epstein birthday schrift}
\pagenumbers{551}{576}
\received{15 November 1997}
\published{29 October 1998}
\papernumber{26}

\newtheorem{theoreme}{{Th\' eor\`eme}}[section]
\newtheorem{theoreme.algebrique}[theoreme]{{Th\' eor\`eme
alg\'ebrique}}
\newtheorem{theoreme.geometrique}[theoreme]{{Th\' eor\`eme
g\'eom\'etrique}}
\newtheorem{proposition}[theoreme]{{Proposition}}
\newtheorem{proposition.fondamentale} [theoreme]{{Proposition
fondamentale}}
\newtheorem{lemme}[theoreme]{{Lemme}}
\newtheorem{corollaire}[theoreme]{{Corollaire}}
\newtheorem{fait}[theoreme]{{Fait}}

\theoremstyle{definition}
\newtheorem{rremarque}[theoreme]{{Remarque}}
\newtheorem{rremarque.terminologique}[theoreme]{{Remarque
terminologique}}

\newtheorem{eexemple}[theoreme]{{Exemple}}

\newtheorem{definition}[theoreme]{{Definition}}

\begin{document}

\title{Sur les espaces-temps homog\`enes}
\asciititle{Sur les espaces-temps homogenes}
\covertitle{Sur les espaces-temps homog\noexpand\`enes}

\author{Abdelghani Zeghib}

\address{Ecole Normale Superieure de Lyon\\
UMPA, UMR 128 CNRS\\
46 Allee d'Italie, 69364 Lyon, FRANCE}

\email{zeghib@umpa.ens-lyon.fr}

\begin{abstract}
Here, we classify Lie groups acting isometrically on compact Lorentz
manifolds, and in particular we describe the geometric structure of
compact homogeneous Lorentz manifolds.
\end{abstract}



\maketitle

\section{Introduction}
Une vari\'et\'e homog\`ene $M$ est par d\'efinition munie d'une
action transitive d'un groupe de Lie $G$, de telle  fa\c con que
$M$ s'identifie \`a un quotient $G/H$ o\`u $H$ est le groupe d'isotropie
(d'un certain point).  Dans
 la suite on supposera toujours
que l'action de $G$ est fid\`ele.

En g\'en\'eral, l'action  de $G$ pr\'eserve une certaine structure
g\'eom\'etrique ``rigide'' \cite{Gro}. Les plus
 belles de ces structures   sont certainement
 les m\'etriques pseudo-riemanniennes. Parmi ces derni\`eres, on distingue
``dans l'ordre'' le cas riemannien et ensuite le cas lorentzien (i.e.\ une
m\'etrique pseudo-riemannienne de type $(1,n-1)$).

Lorsque
$M=G/H$ est une vari\'et\'e riemannienne
homog\`ene compacte,  $G$ est n\'ecessairement
compact (on avait suppos\'e l'action fid\`ele!). Quant \`a $H$, il peut
\^etre n'importe
quel sous-groupe ferm\'e (pas n\'ecessairement discret)
de $G$.

Il n'en est rien, lorsque $M$ est de type lorentzien (et toujours
suppos\'ee compacte).
 Le groupe $G$ peut bien \^etre non-compact, et
 de plus \'etant donn\'e $G$, il
n'est pas \'evident quels sous groupes ferm\'es $H$ peuvent figurer.

Notre but ici est d'essayer de  comprendre,  comme c'est le cas des
m\'etriques riemanniennes, la structure des vari\'et\'es lorentziennes
 homog\`enes, ayant un volume fini.

\subsection{Exemples}
\subsubsection{Cas semi-simple: $SL(2, {\bf R})$} Pour $G$ semi-simple, sa
forme de Killing
 d\'etermine
une m\'etrique pseudo-rieman\-nienne bi-invariante. Ainsi,
elle passe aux
quotients $ G / \Gamma$, pour $\Gamma$ discret, qui seront de
plus munis d'une action \`a gauche isom\'etrique de $G$.
Cette m\'etrique est lorentzienne exactement lorsque $G$ est
localement isomorphe
\`a $SL(2, {\bf R})$.

\subsubsection{Cas r\'esoluble: Groupes de
Heisenberg tordus} \label{Heisenberg.tordu}

La discussion concernant les exemples qui suivent,  se trouve
en grande partie dans
 \cite{Med-Rev}. Il en a \'et\'e \'egalement question dans \cite{Gro} et \cite{Zim}.

L'alg\`ebre de Heisenberg ${\cal HE}_d$ de dimension
$2d+1$ est identifi\'ee en tant qu'espace vectoriel \`a
${\bf R} \bigoplus {\bf C}^d$. Si $Z$ (resp.
$\{e_1,\ldots,e_d \}$ ) est la base canonique de ${\bf R}$
 (resp. ${\bf C}^d$), alors les  crochets non nuls  sont
donn\'es par: $[e_k,ie_k]=Z$. En d'autres termes, si $\omega$ est
la forme symplectique canonique sur ${\bf C}^d$,
$\omega (X,Y) = \langle X,iY{\rangle }_0$, o\`u $\langle \; ,\; {\rangle }_0$ est le produit hermitien
canonique,  alors
 $[X,Y] = \omega(X,Y)Z$.

Consid\'erons l'alg\`ebre r\'esoluble ${\cal HE}^t_d$ (alg\`ebre
de Heisenberg tordue canonique) d\'efinie en ajoutant
un \'el\'ement ext\'erieur  $t$, v\'erifiant $[t,e_k]=ie_k, [t,ie_k]=-e_k$,
pour
$1 \leq k \leq d$ et $[t,Z] =0$.

D\'efinissons sur ${\cal HE}^t_d$, un produit scalaire
$\langle \; ,\; \rangle $, par les lois suivantes: ${\bf C}^d$
 est muni du produit scalaire induit par son produit
hermitien canonique $\langle \; ,\; {\rangle }_0$ et est  orthogonal au
2--plan engendr\'e par $t$ et $Z$. De plus $\langle t,t\rangle =\langle
Z,Z\rangle =0$
et $\langle t,Z\rangle =1$.

Il est remarquable que ceci est un produit lorentzien (en particulier
non d\'eg\'e\-n\'er\'e),  qui est
$Ad ({\cal HE}^t_d)$--invariant! (i.e.\ pour tout g\'en\'erateur
$u$, $ad_u$ est antisym\'etrique au sens de $\langle \; ,\;
\rangle $).

Notons $\tilde{He^t_d}$ le groupe simplement connexe
d\'etermin\'e par
${\cal HE}^t_d$. On remarquera dans la suite qu'il admet
bien des r\'eseaux co-compacts. Comme dans
le cas semi-simple, les vari\'et\'es lorentziennes
 quotients qu'ils d\'eterminent sont donc homog\`enes, et
leurs  groupes d'isom\'etries contiennent des quotients de
$\tilde {He^t_d}$.

En fait, on le constatera au long de ce texte, ce n'est jamais le
 groupe $\tilde{He^t_d}$ qui agit (fid\`element), mais un
quotient, par un r\'eseau de son centre. Pour l'expliciter, notons
 $\tilde{He_d}$ le groupe de Heisenberg simplement
connexe et $He_d$ son quotient par un r\'eseau (isomorphe
\`a ${\bf Z}$) de son centre (ce quotient est unique \`a
conjugaison pr\`es).  Maintenant quotienter $He^t_d$
par un r\'eseau central, revient \`a quotienter $He_d$ par
le groupe engendr\'e par une puissance
enti\`ere  de $\exp (t)$.  On notera $He^t_d$ le quotient
obtenu \`a l'aide du groupe engendr\'e par $\exp(t)$.
 Tous les autres quotients sont des extensions de $He^t_d$
 par des groupes cycliques finis.

En fait on peut d\'efinir ces quotients
   comme produit semi-direct
 du cercle $S^1$ par $He_d$. Le cercle agit
 par
rotation
sur le facteur ${\bf C}^d$ et trivialement sur le centre  ${\bf R}$.
 Le cas de  $He^t_d$ correspond \`a celui o\`u l'action de
$S^1$ est fid\`ele.

Consid\'erons en g\'en\'eral une action par automorphismes de
$S^1$ sur l'alg\`ebre de Heisenberg ${\cal HE}_d$. Soit
$\exp (s2 \pi R)$ le groupe \`a un param\`etre d'automorphismes
 ainsi d\'et\'ermin\'e sur le quotient
de ${\cal HE}_d$ par son centre, identifi\'e \`a
${\bf C}^d$. Il pr\'es\'erve la forme symplectique
 canonique $\omega$ sur ${\bf C}^d$.
 Mais un groupe compact de transformations symplectiques de
${\bf C}^d$ est conjugu\'e \`a un sous-groupe de $U(d)$.
Il s'ensuit que (apr\`es conjugaison) $R$ est une application
${\bf C}$--lin\'eaire diagonale (dans une base orthonorm\'ee)
\`a valeurs propres
 $\lambda_1i,\ldots, \lambda_ki$,
o\`u les $\lambda_j$ sont
des nombres entiers (car $\exp (2 \pi R) =1$).

\begin{definition} [Groupes de Heisenberg tordus] On appelera groupe
de Heis\-enberg tordu tout produit semi-direct du cercle $S^1$
par $He_d$ tel les entiers $\lambda_j$ soient tous non nuls
 et de m\^eme signe (en d'autres termes les produits de
 valeurs propres de $R$ sont
tous non nulles et de m\^eme signe. IL est \'egalement \'equivalent \`a
 dire que l'application ${\bf C}$ lin\'eaire sym\'etrique
$iR$ admet des valeurs propres (r\'eelles) non nulles de m\^eme signe).
\end{definition}

Evidemment pour $d=1$, on n'obtient rien
 d'autre que les extensions cycliques finis de
 $He^t_1$. Ces groupes  peuvent en fait  se d\'efinir autrement  comme
 extensions centrales non triviales  du groupe des d\'eplacements
directs du plan euclidien (appel\'e parfois groupe d'Euclide)
 par le cercle $S^1$.
\begin{rremarque.terminologique}
{La terminologie ci-dessus
 n'est certaine-\break ment pas id\'eale. En effet il existe, au moins
 pour $d=1$,  des
 terminologies concurentes. Par exemple, en physique, un groupe
 de Heisenberg  tordu (pour $d=1$) est dit groupe oscillateur \cite{Str},
  et dans un autre domaine d'int\'er\^et, il s'appele groupe diamant.
Apparemment,
le terme, groupe de Heisenberg tordu, contient plus d'informations
 math\'ematiques. }
\end{rremarque.terminologique}

Une vari\'et\'e d'exemples de vari\'et\'es lorentziennes homog\`enes
compactes s'obtient \`a partir de:
\begin{proposition} \label{metrique.invariante}

{\rm (i)}\qua Un groupe de Heisenberg tordu admet une m\'etrique
 lorentzienne bi-invariante. R\'eciproquement si
une alg\`ebre de Lie obtenu comme produit semi-direct
de $S^1$ par ${\cal HE}_d$, admet une m\'etrique lorentzienne bi-invar\-iante,
alors cette alg\`ebre  est l'alg\`ebre de Lie
 d'un groupe de Heisenberg tordu.

{\rm (ii)}\qua A indice fini pr\`es, il y a \'equivalence
entre les r\'eseaux
d'un groupe de Heisenberg tordu de dimension $2d+2$ et
ceux du sous-groupe $He_d$, ainsi que ceux de
 $\tilde{He_d}$ (le groupe de Heisenberg
simplement connexe de
dimension $2d+1$).

\end{proposition}

\begin{proof}[Preuve]
(i)\qua Cherchons les conditions que doit v\'erifier une telle m\'etrique
 $\langle \; ,\; \rangle $. D'abord la $Ad({\cal HE}_d)$ invariance
 de $\langle \; ,\; \rangle $ restreinte
 \`a ${\cal HE}_d$
 entra\^{\i}ne que cette restriction est positive, \`a noyau
exactement le centre.

 Les conditions de $Ad({\cal HE}_d)$ invariance de $\langle \;
,\; \rangle $ elle m\^eme (i.e.\ non restreinte)
 sont beaucoup plus fortes. En effet, on peut supposer
que $R = ad_t$ pr\'es\'erve
 ${\bf C}^d$ et consid\'erons  $X,Y$ deux \'el\'ements de
${\bf C}^d$.
Ecrivons
la condition d'antisym\'etrie:
 $\langle ad_Xt,Y\rangle +\langle t,ad_XY\rangle =0$. Donc: $\langle
RX,Y\rangle = \langle t,Z\rangle \omega(X,Y)$
(o\`u $Z$ engendre le centre).
 N\'ecessairement, $\langle t,Z\rangle \neq 0$, car sinon $\langle
\; ,\; \rangle $ admettra
 un noyau non trivial contenant $Z$.

 On voit ainsi appara\^{\i}tre la condition
sur les valeurs propres de  $R$ car la restriction de
$\langle \; ,\; \rangle $ \`a ${\bf C}^d$ est d\'efinie positive.
Si elle est satisfaite, on d\'efinira la m\'etrique sur
${\bf C}^d$ par $\langle X,Y\rangle =
\alpha \omega (X, R^{-1}Y)$, o\`u $\alpha = \langle t,Z\rangle $ est une constante
 non nulle assurant que la m\'etrique ainsi obtenue
 est positive (sur ${\bf C}^d$).

On v\'erifie
 alors que $R$ restreinte \`a  ${\bf C}^d$ est antisym\'etrique.
Pour que $R$ (non restreinte)  soit antisym\'etrique, il
suffit que la condition suivante se r\'ealise:
$\langle ad_tt,X\rangle +\langle t,ad_TX\rangle =0$, i.e.\ $\langle t,RX\rangle = 0$ pour tout $X \in
{\bf C}^d$. Il r\'esulte de la bijectivit\'e de $R$
sur ${\bf C}^d$ que $t$ est orthogonal \`a ${\bf C}^d$.  Enfin, on
choisit: $\langle t,t\rangle = \beta$, un nombre r\'eel quelconque. La
m\'etrique est ainsi compl\'etement  d\'efinie, avec deux
param\`etres de choix, $\alpha$ et $\beta$.

(ii)\qua Soit $G$ un  groupe
de  Heisenberg tordu, obtenu comme produit semi-direct de $S^1$ par
  $He_d$. Ainsi  $He_d$  est co-compact dans $G$, en
particulier un r\'eseau de  $He_d$ est aussi un r\'eseau dans $G$.
La
 proposition signifie que r\'eciproquement un r\'eseau
de $G$
 coupe $He_d$ en un r\'eseau et aussi qu'un r\'eseau de
$\tilde{He_d}$ se projette sur un r\'eseau de $He_d$.
Ce sont deux faits
 standard de la th\'eorie
 des goupes discrets dont on peut extraire une preuve
de \cite{Rag} (par exemple le premier fait d\'ecoule d'un \'enonc\'e
g\'en\'eral affirmant  qu'un r\'eseau d'un groupe de Lie
r\'esoluble coupe le nilradical en
un r\'eseau de ce dernier).
\end{proof}

Ainsi, concr\`etement, comme
dans le cas pr\'ec\'edent de
$SL(2, {\bf R})$,  les r\'eseaux des groupes de Heisenberg
(simplement connexes), qu'on comprend parfaitement,  permettent
de construire des vari\'et\'es lorentziennes compactes homog\`enes dont
le groupe
d'isom\'etries est (essentiellement) un groupe de Heisenberg tordu.

Remarquons cependant que si l'on quotiente un
groupe de Heisenberg tordu par un r\'eseau $\Gamma$ contenu
 (pas seulement \`a indice fini pr\`es) dans $He_d$, alors
 on n'aura besoin que de l' $Ad(\Gamma)$--invariance de
$\langle \; ,\; \rangle $. Par Zariski densit\'e des r\'eseaux de $He_d$, ceci
\'equivaut au fait que $\langle \; ,\; \rangle $ est $ad( {\cal HE}_d)$--invariante.

\begin{definition} \label{essentiel}
 On dira qu'une m\'etrique lorentzienne sur l'alg\`ebre de Lie
 d'un groupe de Heisenberg tordu, est essentiellement
bi-invariante, si elle est $ad( {\cal HE}_d)$--invariante.

\end{definition}
\begin{rremarque}{
 D'apr\`es la preuve ci-dessus, une m\'etrique essentiellement
 bi-invariante v\'erifie les m\^emes conditions qu'une   m\'etrique
 bi-invariante, sauf celle de l'othogonalit\'e de $t$ \`a
 ${\bf C}^d$. Une telle
m\'etrique d\'epend donc des deux param\`et\-res, $\alpha$ et  $\beta$, ainsi
 que $2d$ autres param\`etres donnant le produit
de $t$ avec les \'elements d'une (${\bf R}$--) base de
${\bf C}^d$. }
\end{rremarque}

\subsection{Classification}

Notons que malgr\'e  son importance (du moins math\'ematique),
 en dehors des exemples de \cite{Med-Rev} signal\'es ci-dessus,
le seul r\'esultat sensible
connu au  sujet des vari\'et\'es lorentziennes
homog\`enes, est
celui de \cite{Mar},
affirmant que les vari\'et\'es lorentziennes homog\`enes compactes
(ou plus
g\'en\'eralement pseudo-riemannien\-nes) sont g\'eod\'esiquement compl\`etes.
 On peut aussi noter la classification par  \cite{Wol} des vari\'et\'es
lorentziennes
homog\`enes \`a courbure constante, mais pas n\'ecessiarement compactes,
ainsi que le r\'esultat de \cite{D'A} affirmant (entre autres) qu'une vari\'et\'e
lorentzienne homog\`ene compacte et  simplement
connexe est de type riemannien. (On reviendra plus loin au cas
non homog\`ene, o\`u on citera surtout \cite{Zim} et \cite{Gro}).

Le but de
cet article est de montrer que les
 exemples pr\'ec\'edents sont essentiellement les seuls:

\begin{theoreme}   Un espace-temps homog\`ene, de
 volume fini, qui n'est
 pas de type riemannien, admet un sous-groupe normal
 co-compact dans son groupe d'isom\'etries g\'en\'eral,
 qui est soit un rev\^etement fini de
$PSL(2, {\bf R})$, soit
un  groupe de Heisenberg tordu. L'alg\`ebre de Lie
 de ce sous-groupe est en fait un facteur direct dans l'alg\`ebre
de tous les champs de Killing.  De plus ce  sous-groupe
 agit localement librement.
\end{theoreme}

Ce r\'esultat nous permet entre autres de r\'epondre \`a la question
 qu'on s'\'etait pos\'ee pr\'ec\'edemment: si $M =G /H$, alors
quels sous-groupes d'isotropie $H$ peuvent figurer?  Il d\'ecoule
 du th\'eor\`eme pr\'ec\`edent que $H$ est essentiellement discret
au sens que sa composante neutre est compacte.
 Le r\'esultat suivant explicite compl\`etement
la structure topologique
et g\'eom\'etrique des vari\'et\'es lorentziennes
 homog\`enes.

\begin{theoreme} [Classification] \label{classification}
  Soit $M$ une
vari\'et\'e lorentzienne
homog\`ene de volume fini. Supposons que  $M$ n'est
pas  de type riemannien (i.e.\ \`a groupe d'isom\'etries compact). Alors:

{\rm(i)}\qua ou bien $Isom(M)$ contient un rev\^etement
 fini de $PSL(2, {R})$. Dans ce cas
$M$ admet un rev\^etement isom\'etrique
$\tilde{M}$ qui est un produit m\'etrique
 de $\widetilde{ SL(2, {\bf R}) }$ muni de sa forme
de Killing, par une vari\'et\'e riemannienne
homog\`ene compacte.

{\rm(ii)}\qua ou bien $Isom(M)$ contient $S$ un groupe
 de Heisenberg tordu.
  Dans ce cas $M$ admet un rev\^etement $\tilde{M}$
qui  se construit
 de la fa\c con suivante. Il existe
  une vari\'et\'e riemannienne homog\`ene
 compacte $(L, r)$,  munie d'une
action isom\'etrique localement libre de $S^1$.
 Le cercle $S^1$ isomorphe au centre de $S$,
y agit par translation et agit par suite diagonalement
sur $S \times L$, muni de la  m\'etrique produit
de celle de L par une m\'etrique lorentzienne essentiellement bi-invariante
sur $S$.
Alors le rev\^etement $\tilde{M}$ est le quotient
$ S \times_{S^1} L$  de cette action. Il est
muni  de la m\'etrique
d\'eduite par projection,  de la m\'etrique induite sur
$ TS \oplus {\cal O}$,  o\`u
 $\cal O$ est la distribution orthogonale
\`a l'action de $S^1$ sur $L$.

En fait $M= \tilde{M} / \Gamma$, o\`u
$\Gamma$ est le graphe d'un homomorphisme $\rho$
d'un r\'eseau co-compact $\Gamma_0$ de $S$
 dans  $Isom_{S^1}(L)$,
 le groupe
d'isom\'etries de $L$ respectant l'action de $S^1$. De plus
 le centralisateur de $\rho (\Gamma_0)$ dans
$Isom_{S^1}(L)$ agit transitivement
 sur $L$.

\end{theoreme}

On peut par exemple prendre pour $L$ la sph\`ere $S^3$ munie
d'une fibration de
Hopf. Le groupe
d'isom\'etries qui la pr\'eserve
est isomorphe  \`a $S^1 \times S^3$.
 On prendra pour  $\rho$ un homomorphisme  d'un r\'eseau
de $S$ \`a valeurs dans $S^1$ (ce qui
assurera que le centralisateur de l'image
de $\rho$ agit transitivement sur $S^3$). Le groupe d'isom\'etries
de la vari\'et\'e lorentzienne homog\`ene
 compacte  ainsi obtenue,
 sera essentiellement
 $S^3 \times S^1 \times S$.

\begin{rremarque} {
On d\'eduit du th\'eor\`eme de   structure ci-dessus qu'on
peut changer la m\'etrique tout en la gardant homog\`ene, de telle
fa\c con que la m\'etrique sur $S$ soit bi-invariante.

 Il est bien connu que sur $SL(2, {\bf R})$,
la forme de Killing est \`a une constante pr\`es, la seule
m\'etrique bi-invariante.  Sur un groupe
de Heisenberg tordu $S$,  il y en a beaucoup, mais
 elles sont toutes isom\'etriques (pas seulement conformes!)
par
automorphismes dans le rev\^etement
universel $\tilde {S}$. Ceci est li\'e au fait (remarquable) qu'une
structure Lorentzienne bi-invariante donn\'ee sur
$\tilde{S}$,  admet des transformations conformes non triviales.
 Elles sont en fait des homoth\'eties, i.e.\ \`a distorsion
 constante.

Ainsi  sur $S$ toutes
 les m\'etriques bi-invariantes  sont localement
isom\'etriques. Cependant il y a un module de dimension 2
 (les
 param\`etres $\alpha$ et $\beta$ de la preuve
pr\'ec\'edente)
 de telles structures globales (voir \ref{bi-invariante}). }
\end{rremarque}

\subsection{Ingr\'edients  de la preuve}

La finitude du volume sera utilis\'ee, comme d'habitude,  pour en d\'eduire
des propri\'et\'es de r\'ecurrence de l'action de $G$. Mais le plus
grand int\'er\^et de
cette hypoth\`ese pour nous ici, c'est de permettre de construire
 un produit scalaire $L^2$, sur l'alg\`ebre de
 Lie de $G$, ayant la propri\'et\'e  \'el\'ementaire mais remarquable
 d'\^etre  bi-invariant.

 En effet, plus g\'en\'eralement, si $G$ est un groupe de
 Lie agissant
sur une vari\'et\'e $M$ (pas
 n\'ecessairement transitivement) en pr\'eservant une m\'etrique
 pseudo-riemannienne $\langle \; ,\; \rangle $,
alors la forme $\kappa(X,Y) = \int_M \langle X(x),Y(x)\rangle dx$ d\'etermine une forme
bilin\'eaire sur l'alg\`ebre de Lie $\cal G$, qui est bi-invariante. Pour
le voir, il
suffit de remarquer que  si $\phi^t$ est un groupe \`a param\`etre de $G$
(identifi\'e au flot correspondant de $M$)
et $X$ est un   \'el\'ement $\cal G$ (identifi\'e au champ de vecteurs
correspondant sur $M$), alors $\phi^t_*X= Ad(\phi^t)X$.

Notons qu'il est aussi  possible de consid\'erer des formes du type
$\kappa(X,Y) = \int_U \langle X(x),Y(x)\rangle dx$, o\`u $U \subset M$ est
 un sous-ensemble
$G$--invariant quelconque, ou plus g\'en\'eralement en int\'egrant
par rapport \`a une mesure $G$--invariante quelconque.
Remarquons aussi  que la m\^eme construction marche lorsque $G$ pr\'eserve
 un tenseur quelconque sur $M$, et permet ainsi de construire un tenseur
 ``analogue''  bi-invariant sur $\cal G$.

Cependant, le r\'esultat obtenu
est g\'en\'eralement trivial (m\^eme nul!). Ainsi, lorsque $G$ est
simple, la forme obtenue est un multiple (souvent nul) de la forme de
Killing. On peut par exemple  prendre $M=G$, qu'on munit d'une structure
pseudo-riemannienne
invariante \`a gauche (elle s'obtient
simplement d'un produit scalaire de m\^eme signature sur $\cal G$). Ainsi
 $G$ agit sur $M$ en respectant cette structure.  Lorsque
$G$ est compact la forme $\kappa$ construite sur $\cal G$ sera
 d\'efinie positive, d\'efinie n\'egative ou nulle, quelle que
 soit la signature de la structure pseudo-riemannienne de d\'epart.

\subsection{Cas lorentzien}
Dans notre cas lorentzien, la forme $\kappa$ sera suffisamment
non triviale
d\`es qu'il existe des
champs $X \in G$ tels que $\langle X(x),X(x)\rangle $ garde un signe constant. Il se
trouve, comme
cela \'etait \'etabli dans \cite{Ze1} que c'est effectivement le cas
pour tout champ $X$ engendrant un groupe \`a param\`etre $\phi^t$ non
pr\'ecompact,
i.e.\ la fermeture dans $G$ de $\{\phi^t / t \in {\bf R} \}$ n'est pas compact.
  C'est \`a ce niveau l\`a qu'on utilise l'aspect dynamique
 de la finitude du volume.
En fait ona:
\begin{proposition.fondamentale}\label{fondamentale}{\rm\cite{Ze1}}\qua Soit
$(M,\langle \; ,\; \rangle )$ une vari\'et\'e lorentzienne de volume fini.
 Soit $X$ un champ de Killing sur $M$, d\'eterminant un
flot non pr\'ecompact,
alors: $\langle X(x),X(x)\rangle \ge 0$ pour tout $x \in M$. On dira que $X$ est
(partout) non temporel.
\end{proposition.fondamentale}
On en d\'eduit ce fait, qui n'entra\^ \i ne pas tout \`a
fait  que $\kappa$
est lorentzienne, exactement  comme $\langle \; ,\; \rangle $, mais en borne la
d\'eg\'en\'erescence:
\begin{proposition}\label{condition*}{\rm Condition ($*$)}\qua Soit
$\cal P$ un sous-espace vectoriel
de champs de Killing tel pour l'ensemble des \'el\'ements
 $X \in {\cal P}$
d\'eterminant des flots non pr\'ecompacts, est dense dans
$\cal P$. Alors la forme $\kappa$ est positive
sur $\cal P$ et son noyau est au plus de dimension 1 (ou
en d'autres termes, l'ensemble des vecteurs isotropes de
$\cal P$ est un sous-espace vectoriel de dimension $\leq 1$).
\end{proposition}

\subsection{Un r\'esultat alg\'ebrique}

Il se trouve que les donn\'ees alg\'ebriques, fournies par
 $\kappa$, v\'erifiant la proposition pr\'ec\'edente, suffisent
largement  pour
 comprendre le groupe $G$:

\begin{theoreme.algebrique} \label{theoreme.algebrique}
Soit $G$ un groupe de Lie non compact  dont l'alg\`ebre de Lie
$\cal G$ est munie d'une forme bi-invariante $\kappa$,
v\'erifiant l'hypoth\`ese de non d\'eg\'en\'erescence  {\rm($*$)} suivante:

{\rm Condition ($*$)}\qua
Si  $\cal P$ est un sous-espace vectoriel
de $\cal G$,  tel que l'ensemble des \'el\'ements  $X \in {\cal P}$
d\'eterminant des groupes \`a param\`etre
 non pr\'ecompacts est dense dans
$\cal P$; alors la forme $\kappa$ est positive
sur $\cal P$ et son noyau est au plus de dimension 1.

Alors $\cal G$ s'\'ecrit comme  somme directe orthogonale
  d'alg\`ebres:
${\cal G}= {\cal K}+{\cal A}+{\cal S}$, o\`u:
 $\cal K$ est une alg\`ebre compacte (i.e.\ l'alg\`ebre
de Lie d'un groupe de Lie semi-simple compact), $\cal A$ est
une alg\`ebre ab\'elienne, et $\cal S$ est soit triviale, soit
$sl(2, {\bf R})$, soit l'alg\`ebre de Lie du groupe affine (des
transformations de la droite),
soit une alg\`ebre de Heisenberg
 ${\cal HE}_d$, soit une
alg\`ebre de Heisenberg tordue. On a:

{\rm (i)}\qua Lorsque $\cal S$ est  triviale,
 $\kappa$ est positive \`a noyau
de dimension $\leq 1$. Lorsque
  $\cal S$ est non triviale,
 $\kappa$ est d\'efinie positive sur $\cal K$ et $\cal A$.

{\rm (i)}\qua Lorsque $\cal S$ est
l'alg\`ebre
de Lie du groupe affine, $\kappa$ est positive
 d\'eg\'en\'er\'ee sur $\cal S$
 et admet pour  noyau exactement
 l'id\'eal d\'etermin\'e par les translations.

{\rm (i)}\qua Lorsque
$\cal S$ est une alg\`ebre de Heisenberg,  $\kappa$ est positive
 d\'eg\'en\'er\'ee
sur $\cal S$ et admet pour  noyau exactement le centre.

{\rm(iv)}\qua Lorsque  $\cal S$ est  une alg\`ebre de Heisenberg
tordue, la forme $\kappa$ est lorentzienne sur $\cal S$. Le
sous-groupe de $G$ d\'etermin\'e par $\cal S$ est un groupe de
 Heisenberg tordu.  De plus le sous-groupe abelien d\'etermin\'e
 par ${\cal A}+{\cal Z}$, o\`u $\cal Z$ est le centre
de $\cal S$, est compact.

{\rm(v)}\qua Lorsque ${\cal S} = sl(2, {\bf R})$, la forme $\kappa$
 sur $\cal S$ est lorentzienne et le sous-groupe d\'etermin\'e
 par $\cal S$ est
un rev\^etement
fini de $PSL(2, {\bf R})$. De plus le sous-groupe d\'etermin\'e par
 $\cal A$ est compact.

\end{theoreme.algebrique}

\subsection{Un r\'esultat g\'eom\'etrique}

Le th\'eor\`eme alg\'ebrique s'applique
en particulier
aux  groupes de Lie connexes non compact
 agissant isom\'etriquement sur une vari\'et\'e Lorentzienne
 $(M,\langle \; ,\; \rangle )$ de volume fini. Certaines parties de ce
th\'eor\`eme sont dues dans ce
cas \`a \cite{Zim} et ensuite \cite{Gro}. Plus pr\'ecis\'ement, la structure
alg\'ebrique de $\cal G$ est \'elucid\'ee dans \cite{Zim} lorsque
 $\cal G$ contient $SL(2, {\bf R})$. Il y a \'et\'e \'egalement
d\'emontr\'e que le nilradical est de degr\'e de nilpotence $\leq 2$.

Dans \cite{Gro}, il a \'et\'e question d'am\'eliorations g\'eom\'etriques
 de r\'esultats de \cite{Zim} (surtout dans le cas analytique).
En effet, on peut, en g\'en\'eral,  am\'eliorer le th\'eor\`eme
alg\'ebrique pr\'ec\'edent, par un r\'esultat g\'eom\'etrique, ponctuel.
Il exprime
 essentiellement le fait que si un champ de Killing $X$ est non temporel
au sens de $\kappa$ (i.e.\ $\kappa (X,X) \ge 0$), c'est qu'il l'est
 ponctuellement au sens
 de $\langle \; ,\; \rangle $ (i.e.\ $\langle X(x),X(x)\rangle   \geq 0$, pour tout $x \in M$).

Tout ce
qui concerne $SL(2, {\bf R})$ dans les r\'esultats suivants est
d\'emontr\'e par
\cite{Gro}.  Notre approche ici permet de les red\'emontrer.

\begin {theoreme.geometrique}
Soit $G$ un groupe de Lie connexe non compact
agissant isom\'etriquement
sur une vari\'et\'e lorentzienne $(M,\langle \; ,\; \rangle )$ de volume fini. Notons
 $\kappa$ la forme ainsi d\'efinie sur $\cal G$.

{\rm1)}\qua Supposons que $\kappa$ est positive, alors les orbites de $G$ sont
 non temporelles (i.e.\ la restriction de $\langle \; ,\; \rangle $ \`a ces orbites est $\ge 0$).
Le noyau de $\kappa$, s'il n'est pas trivial est un champ de Killing
(partout) de type lumi\`ere (au sens de
$\langle \; ,\; \rangle $) et \`a orbites g\'eod\'esiques.

{\rm2)}\qua Supposons que $\kappa$ n'est pas positive, donc $\cal G$ contient un
facteur direct $\cal S$, isomorphe \`a $sl(2, {\bf R})$ ou alg\`ebre de
 Heisenberg tordue. Alors l'action de $\cal S$ est
partout localement libre.
\end{theoreme.geometrique}

Le r\'esultat de \cite{Gro}  pour $sl(2, {\bf R})$  est
plus pr\'ecis.  Il
affirme davantage que la distribution orthogonale (aux orbites)
 est int\'egrable
(et aussi g\'eod\'esique). Il s'ensuit qu' un certain rev\^etement est
un produit tordu d'une vari\'et\'e riemannienne par
 $\widetilde{SL(2, {\bf R})}$.

En fait lorsqu'un groupe isomorphe \`a   $\widetilde{SL(2, {\bf R})}$
 ou \`a un groupe
 de Heisenberg tordu agit isom\'etriquement sur
une vari\'et\'e lorentzienne de volume fini, alors celle ci s'obtient
pratiquement de la m\^eme fa\c con que dans le cas homog\`ene,
explicit\'e au th\'eor\`eme \ref{classification}, \`a ceci pr\`es que
$L$ ne sera suppos\'ee ni homog\`ene ni compacte:

\begin{theoreme}\label{structure.semisimple}{\rm\cite{Gro}}\qua Une vari\'et\'e
lorentzienne de volume fini
 munie d'une action isom\'etrique d'un
groupe localement isomorphe \`a $SL(2, {\bf R})$ est rev\^etue par un
produit de $ \widetilde {
 SL(2, {\bf R}) }$ par une vari\'et\'e riemannienne  $ (L, r)$,
muni d'une m\'etrique tordue $h_{(g,x)} = f(x) k \otimes r_x$,
 o\`u $f$ est une fonction positive sur $L$ et $k$ est la
forme de Killing de $\widetilde {SL(2, {\bf R})}$.

\end{theoreme}

Ici on a  un r\'esultat de structure,
un peu plus compliqu\'e, dans le cas d'un groupe de
Heisenberg tordu $G$
, du fait que la distribution
orthogonale
n'est pas n\'ecessairement int\'egrable. C'est en fait
 sa satur\'ee
 par le centre de $G$ qui l'est.

La construction est la suivante.
 Soit $(L, r)$ une vari\'et\'e riemannienne munie d'une
action isom\'etrique localement libre de $S^1$. Notons
 $\cal
 O$ la distribution orthogonale aux orbites.

Soit $\cal M$ l'espace des m\'etriques lorentziennes
 essentiellement  bi-invariantes sur ${\cal G}$ (\ref{essentiel}).
 Soit $\phi\co L \to {\cal M}$ une application (de m\^eme
 classe de r\'egularit\'e que toutes
les donn\'ees). Munissons le
produit $ G \times L$ de la m\'etrique tordue d\'efinie par
$\phi$:  $h_{(g,x)} = m_x \otimes r_x$ (l'espace tangent au facteur
$G$ \'etant partout identifi\'e \`a $\cal G$).

Le centre de $G$, isomorphe \`a $S^1$,
 y agit isom\'etriquement par translation.
On a donc une action isom\'etrique
 diagonale de $S^1$ sur le produit
 $ G \times L$. Notons
 $G  \times_{S^1} L$ le quotient et
 munissons le de la m\'etrique
d\'eduite par projection,  de la m\'etrique induite sur
 l'horisontal $ {\cal G} \oplus {\cal O}$

Soit $\Gamma$ un r\'eseau de $G \times Isom_{S^1}(L)$
  o\`u $Isom_{S^1}(L)$ d\'esigne
 le groupe d'isom\'etries de $L$ pr\'eservant
l'action de $S^1$. On suppose que $\Gamma$ agit
 sans point fixe sur $G \times _{S^1} L$ ( ce
qui sera toujours vrai pour un sous-groupe
d'indice fini). Le quotient $M= \Gamma \setminus G \times_{S^1} L$
 est  une  vari\'et\'e lorentzienne de
volume fini munie d'une action isom\'etrique de  $G$.

\begin{theoreme} \label{structure.resoluble} Toute vari\'et\'e lorentzienne
de volume
fini, munie d'une action isom\'etrique d'un groupe de Heisenberg
tordu $G$ est
construite de la fa\c con pr\'ec\'edente.
\end{theoreme}

\begin{eexemple} { On peut prendre pour $L$ le groupe
$G$ lui m\^eme muni d'une m\'etrique
riemannienne invariante \`a
droite, et de l'action de son centre. On voit
 sur  cet exemple que $\cal O$ peut bien \^etre non int\'egrable.
 En jouant sur $\Gamma$, qui est un r\'eseau
 de $G \times G$, on peut r\'ealiser diverses
propri\'et\'es de densit\'e des orbites de $G$. }
\end{eexemple}

{\it
 La classification des alg\`ebres de Lie de groupes agissant
isom\'etriquement sur des vari\'et\'es lorentziennes compactes, a \'et\'e
d\'emontr\'ee ind\'ependemment par S Adams et G Stuck \cite{AS1}. Le
 pr\'esent article, ainsi que \cite{AS1} sont parus simultan\'ement
(sous forme de preprints) en Mai 1995. D'autres r\'esultats
compl\'ementaires qui pr\'ecisent cette classification
ont \'et\'e ensuite d\'emontr\'es dans \cite{AS2} et
\cite{Ze2}.}

\section{La condition ($*$)}

Rappelons bri\`evement dans ce qui suit  les \'el\'ements de la  preuve de
 \ref{fondamentale} \cite{Ze1}. Le premier point  est que dans
 un groupe de Lie la fermeture $L$ d'un groupe \`a un param\`etre
 $\phi^t$ est soit ${\bf R}$ soit un tore (compact). En
effet $L$ est un groupe ab\'elien poss\'edant un groupe \`a un
param\`etre dense.

Il en d\'ecoule que si une sous-suite $\phi^{t_i}$ est
pr\'ecompacte (i.e.\ \'equicontinue) alors le flot
$\phi^t$ lui m\^eme est pr\'ecompact.

Le second point est un ph\'enom\`ene d'uniformit\'e
valable pour des suites de transformations $f_i$
 pr\'eservant une connexion. Il stipule que
 si la suite des d\'eriv\'es $D_xf_i$ en un point $x$
 donn\'e est \'equicontinue, alors la suite $f_i$ elle m\^eme
est \'equicontinue. Ceci d\'ecoule de la d\'efinition m\^eme
de la structure diff\'erentiable du groupe $G$
d'isom\'etries de la connexion. En effet cette structure est
 caract\'eris\'ee par le fait que pour tout rep\`ere $r_x$ en $x$,
 l'\'evaluation $e\co G \to  Rep(M)$, $e(f) = f^*(r_x)$
 est un plongement propre.

Le dernier point est qu'au voisinage  d'un point $x$,
 qu'on
peut supposer r\'ecurrent, o\`u
le champ de Killing $X$ g\'en\'erateur
infinit\'esimal de $\phi^t$ est de type temps,
les applications  de retour, ont leurs d\'eriv\'ees \'equicontinues
en $x$. En effet, ces d\'eriv\'ees respectent la m\'etrique
 riemannienne (d\'efinie au voisinage de $x$) obtenue
 canoniquement \`a partir de la m\'etrique lorentzienne,
juste en changeant le signe le long de $X$.  \endproof \\

La Proposition \ref{condition*} d\'ecoule du fait suivant:

\begin{lemme} Soit $\cal P$ un sous-espace vectoriel
de champs de Killing tel que pour tout $X \in {\cal P}$ et $x \in M$,
$\langle X(x),X(x)\rangle  \ge 0$. Alors la forme $\kappa$ est positive
sur $\cal P$ et son noyau est au plus de dimension 1.
\end{lemme}
\begin{proof}[Preuve]
Il d\'ecoule de
l'hypoth\`ese que si $X \in {\cal P}$ est isotrope au sens
de $\kappa$, alors $X(x)$ est isotrope au sens de $\langle \;
,\; \rangle _x$
pour tout $x$. Donc si ${\cal A}$ est un sous-espace isotrope
 de ${\cal P}$, alors:
${\cal A}_x = \{X(x), X \in {\cal A} \}$ est un sous-espace isotrope
de $(T_xM, \langle \; ,\; {\rangle }_x)$. Il s'ensuit que: dim$({\cal A}_x) \leq 1$
 pour tout $x$ car la m\'etrique $\langle \; ,\; \rangle $ est lorentzienne.

La preuve du  lemme
sera achev\'ee si l'on montre que deux champs de Killing (partout)
colin\'eaires, sont tels que l'un est multiple constant de l'autre.

En effet soit $X$ et $Y$
 deux tels champs et  \'ecrivons (localement) $Y =fX$ o\`u $f$ est une
certaine fonction.
Notons $\nabla$ la connection de Levi-Civita. Alors, un Champ
de Killing tel que $X$ v\'erifie que: pour tout $x$, l'application
$u \in T_xM \to \nabla_uX \in T_xM$ est antisym\'etrique.
Ainsi $0= \langle \nabla_u(fX), u\rangle = (u.f)\langle X,u\rangle +f\langle \nabla_uX,u\rangle =(u.f)\langle X,u\rangle $, car
$X$ et $Y=fX$ sont, tous les deux, deschamps de Killing. Il en d\'ecoule que $f$
est constante.
\end{proof}

\section{D\'ebut de la preuve du th\'eor\`eme alg\'ebrique}

Sans le mentionner,
on utilisera parfois, l'affirmation suivante, qui contient des
 faits classiques standards:

\begin{fait}  \label{folk} Soit ${\cal G}$ une alg\`ebre de
Lie muni d'une forme bi-invariante $k$. On a:

{\rm(i)}\qua Le noyau de $k$ est un id\'eal de ${\cal G}$.

{\rm(ii)}\qua Si $k$ est d\'efinie positive, alors
${\cal G}$ est somme directe $k$--orthogonale
d'une alg\`ebre abelienne et d'une alg\`ebre
compacte (i.e.\ l'alg\`ebre de Lie
d'un groupe de Lie semi-simple compact).

{\rm(iii)}\qua Si ${\cal G} $
est compacte, alors $k$ est multiple de sa forme de Killing.

\end{fait}

Soit maintenant  ${\cal G}$ une alg\`ebre de Lie
munie d'une forme $\kappa$ comme  dans le the\'eor\`eme alg\'ebrique.

\begin{lemme} \label{dense.abelien} Soit $\cal P$ une sous alg\`ebre
ab\'elienne
de ${\cal G}$
 ayant un \'el\'ement $X$
d\'eterminant un
flot non pr\'ecompact.  Alors la forme $\kappa$ est positive
sur $\cal P$ et son noyau est au plus de dimension 1.
\end{lemme}
\begin{proof}[Preuve]
On applique la condition ($*$) sachant que la fermeture
du groupe d\'etermin\'e par $\cal P$ est un produit d'un tore par un
espace vectoriel non trivial. Tous les groupes \`a un param\`etre
sont non pr\'ecompacts sauf ceux tangents au facteur torique.
\end{proof}

\paragraph{Le nilradical}

Le lemme \ref{dense.abelien} s'\'etend en fait aux groupes
 nilpotents gr\^ace \`a la:

\begin{proposition} \label{dense.nilpotent} L'ensemble des groupes
\`a un param\`etre non pr\'ecompacts d'un groupe de Lie nilpotent
non compact, est dense. C'est en fait le compl\'emen\-taire d'un
tore maximal (qui est par ailleurs central et unique).
\end{proposition}

\begin{proof}[Preuve] Soit $N$ un tel groupe,
  $\tilde{N}$  son groupe rev\^etement
universel,  et $\Gamma$  le groupe fondamental de $N$. Alors
$\Gamma$ est
central dans $\tilde{N}$. De plus, c'est un r\'eseau dans
 un unique
 sous-groupe
 de Lie (connexe) $\tilde{L}$, \'egalement central
(pour d\'efinir
$\tilde{L}$, on se ram\`ene au cas ab\'elien,
en   remarquant simplement  que le centre de $\tilde{N}$ est
connexe, car si un \'el\'ement est central, alors  le
groupe \`a param\`etre
(unique) qui le contient est central). La projection
de $\tilde{L}$ dans $N$ est un tore (maximal).

Ainsi $N$ se projette sur $N/L = \tilde{N} /\tilde{L}$, qui est simplement
 connexe,
donc ayant tous ses groupes \`a un param\`etre non pr\'ecompacts. Il
 s'ensuit que
les groupes \`a un param\`etre de $N$ qui sont pr\'ecompacts, sont
tangents \`a l'alg\`ebre de Lie de $L$.
\end{proof}

Notons $N$ le  nilradical de $G$, i.e.\ le plus
grand sous groupe de Lie (connexe) normal  nilpotent. On
supposera dans cette section qu'il est non compact.

\begin{corollaire} \label{nilradical} Si le nilradical $N$ est non
compact, alors
 la restriction de $\kappa$ \`a $\cal N$ est une forme
 positive, dont le noyau est
un id\'eal $\cal I$
de dimension
 $\leq 1$. De plus $\cal N$ est isomorphe \`a
 une somme directe orthogonale d'alg\`ebres ${\cal N}=
{\cal A}
 + {\cal HE}_d$,  o\`u $\cal A$ est ab\'elienne et
 ${\cal HE}_d$ est l'alg\`ebre de Heisenberg
de dimension $2d+1$.

L'action adjointe de $G$ sur
${\cal N}/{\cal I}$ est  \`a image compacte (car elle pr\'eserve
une forme d\'efinie positive).

Lorsque le facteur  correspondant \`a
l'alg\`ebre de Heisenberg est non trivial, le noyau $\cal I$ de
$\kappa$ est exactement son centre ${\cal Z}$.
\end{corollaire}

\begin{proof}[Preuve]
On utilise \ref{dense.abelien} pour en d\'eduire que
$\kappa$ est positive sur ${\cal N}$ et
que
dim${\cal I} \leq 1$.
On remarque ensuite  que l'alg\`ebre
${\cal N} / {\cal I}$ est ab\'elienne, car elle est nilpotente et
admet une m\'etrique d\'efinie positive bi-invariante.
\end{proof}

\begin{rremarque} \label{non.canonique} { $\cal A$ n'est pas
canoniquement d\'efinie,
mais la somme
 ${\cal A}+{\cal Z}$ et le facteur de type Heisenberg
${\cal HE}_d$  le sont. }
\end{rremarque}

\begin{proposition} \label{radical}

{\rm(i)}\qua Le centre $\cal Z$ de ${\cal HE}_d \subset {\cal N}$ est en fait
central dans $\cal G$.

{\rm(ii)}\qua Tout $X \in {\cal HE}_d \subset {\cal N}$ non central, engendre un groupe
\`a un param\`etre non pr\'ecompact.

{\rm(iii)}\qua Si $Y$ est un \'el\'ement non trivial de
$\cal G$ qui commute avec un \'el\'ement non
central de $\cal N$, alors $\kappa(Y,Y) > 0$.
\end{proposition}

\begin{proof}[Preuve]
(i)\qua Soit $A$ un automorphisme de ${\cal HE}_d$ respectant
 $\kappa$. Supposons que $A$ induit sur $\cal Z$ une homoth\'etie
non triviale.
  Alors $\cal Z$ sera le seul sous espace propre associ\'e
\`a une valeur propre de module $ \neq 1$, car
$\kappa$ est d\'efinie positive sur ${\cal HE}_d /{\cal Z}$.
 Il existera donc un suppl\'ementaire $\cal T$ de $\cal Z$, sur lequel
$A$ respecte une
m\'etrique d\'efinie positive (et en particulier \`a
 valeurs propres de module \'egale \`a 1). On obtient une contradiction
en consid\'erant deux \'el\'ements,
$X$ et $Y$ de $\cal T$, v\'erifiant $[X,Y] = Z \in {\cal Z}$.

(ii)\qua d\'ecoule du fait qu'alors $ad_X$ est nilpotent (et
 non trivial) et
donc le groupe \`a un
param\`etre $exp(t ad_X)$ est non pr\'ecompact.

(iii)\qua En effet si $Y$ commute avec un \'el\'ement non central
$X \in {\cal HE}_d$, alors $Y$, $X$ et $\cal Z$ d\'eterminent
une sous-alg\`ebre ab\'elienne
de dimension $\geq 2$, v\'erifiant \ref{dense.abelien}. Il s'ensuit que $\cal Z$
 est le seul espace $\kappa$ isotrope de cette sous-alg\`ebre.
\end{proof}

\begin{proposition} \label{semisimple.nilpotent}
Soit ${\cal L} \subset {\cal G}$ une sous-alg\`ebre
semi-simple. Alors la somme ${\cal G}^\prime= {\cal L}+{\cal N}$ est
orthogonale (au sens de $\kappa$) et
directe (au sens d'alg\`ebres)
\end{proposition}

\begin{proof}[Preuve]
Soit ${\cal I} \subset {\cal N}$
le noyau de la restriction de $\kappa$ \`a $\cal N$.
 C'est un id\'eal de ${\cal G}^\prime$ de dimension $\leq 1$.
Par semi-simplicit\'e, $\cal L$ centralise $\cal I$. On peut appliquer le
Fait \ref{orthogonal} pour voir que $\cal I$ est orthogonale \`a $\cal L$
 et par suite \`a ${\cal G}^\prime$.

On peut donc en passant au quotient
${\cal G}^\prime /{\cal I} = {\cal L} + {\cal N} / {\cal I}$,
supposer que la forme $\kappa$ est d\'efinie positive
sur ${\cal N}$.

La proposition est bien connue lorsque $\kappa$ est
d\'efinie
positive sur ${\cal G}^{\prime}$ (voir \ref{folk}. On va essayer donc de se
ramener \`a cette situation.
Par
\ref{nilradical}, l'action adjointe de $\cal L$ sur ${\cal N}$
est \`a image compacte. On peut donc supposer que $\cal L$ est
compacte. De plus, quitte \`a tra\^ \i ter facteur par
facteur, on peut supposer que
$\cal L$ est simple. Soit $k$ la forme
de Killing de ${\cal G}^{\prime}$. Elle est
triviale sur ${\cal N}$ car c'est le nilradical
et sa restriction \`a $\cal L$ est
un multiple
non nul de la forme de Killing de $\cal L$. Ainsi,
sur $\cal L$, $\kappa$ est multiple de
$k$. Il s'ensuit qu'il existe un
choix d'un r\'eel
$\alpha$ tel que $\kappa + \alpha k$ soit
d\'efinie positive sur ${\cal G}^{\prime}$. Ainsi, on s'est ramen\'e
au
cas o\`u
$\kappa$
est d\'efinie positive sur ${\cal G}^\prime$.

Pour montrer l'orthogonalit\'e de la somme ${\cal L} +{\cal N}$,
  on utilise le fait g\'en\'eral suivant,
dont la preuve d\'ecoule de la bi-invariance de $\kappa$.
\end{proof}

\begin{fait} \label{orthogonal} Soit $\cal L$ une sous-alg\`ebre de $\cal G$,
et $Y$ un \'el\'ement de $\cal G$ centralisant $\cal L$.
Alors l'application $X \in {\cal L} \to \kappa(X,Y) \in
 {\bf R}$ est un homomorphisme, i.e.\ $\kappa([X,X^\prime],Y)=0$, pour $X$, $X^\prime$ dans
$\cal L$.
\end{fait}

\paragraph{Le radical}

Soit $R$
 le radical   de $G$ (i.e.\ le plus
grand sous groupe de Lie normal r\'esoluble) et $\cal R$
 son alg\`ebre de Lie. On
supposera dans cette section qu'il est non
 compact. On a d'abord la constatation suivante:
\begin{fait}
Si $R$ est non compact, alors le nilradical $N$
l'est \'egalement.
\end{fait}

\begin{proof}[Preuve] En effet,
 s' il est compact, $N$ sera  central dans $G$ et en particulier
dans $R$. Ainsi l'avant dernier groupe
d\'eriv\'e de $R$ contient strictement $N$ et est nilpotent. Par
naturalit\'e, il est
normal dans $G$, ce qui contredit la d\'efinition de $N$.
\end{proof}

La proposition \ref{semisimple.nilpotent} se g\'en\'eralise \`a $\cal R$:

\begin{proposition} \label{semisimple.radical} Soit ${\cal L} \subset {\cal
G}$ une sous-alg\`ebre
semi-simple. Alors la somme ${\cal G}^\prime= {\cal L}+{\cal R}$ est
orthogonale (au sens de $\kappa$) et
directe (au sens d'alg\`ebres)
\end{proposition}

\begin{proof}[Preuve] Cela d\'ecoule de \ref{semisimple.nilpotent}
 et du lemme suivant. \end{proof}
\begin{lemme}
Soit $A$ un automorphisme semi-simple de $R$, trivial sur $\cal N$,
 alors $A$ est trivial.
\end{lemme}
\begin{proof}[Preuve]
 Soit ${\cal E} \subset {\cal R}$ un sous-espace
vectoriel suppl\'ementaire de $\cal N$ invariant par $A$.
Soit $X \in {\cal E}$, $Y \in {\cal N}$, alors
$[X,Y] \in {\cal N}$. Donc $[X,Y] = A[X,Y] = [AX,Y]$. Autrement
dit $X-AX$ centralise $\cal N$. Par maximalit\'e de $\cal N$
en tant que sous-alg\`ebre normale nilpotente, on d\'eduit
que l'application $X \in {\cal E} \to X-AX \in {\cal E}$ est
nulle (car son image est contenue dans $\cal E$).
\end{proof}

\paragraph{Facteur semi-simple}

\begin{fait} \label{facteur.semisimple} Supposons que $R$ est non compact.
Alors on a une d\'ecomposition directe et orthogonale
 ${\cal G} = {\cal K} +{\cal R}$ o\`u
$\cal K$ est une sous-alg\`ebre semi-simple compacte.
\end{fait}
\begin{proof}[Preuve]
 D'apr\`es ce qui pr\'ec\`ede, il suffit simplement de
montrer que le facteur semi simple  $\cal K$ est compacte. Il suffit
pour cela d'observer que la restriction de $\kappa$
\`a chaque facteur de $\cal K$ est positive et non
triviale.  Pour cela
 on applique la condition ($*$) \`a l'alg\`ebre
 ${\cal K}^\prime={\cal K} + {\bf R}X$, o\`u $X$ est un \'el\'ement de
 $\cal R$  qui d\'etermine un groupe \`a un param\`etre
 non pr\'ecompact. En effet tous les groupes \`a un param\`etre
 de ${\cal K}^\prime$ non tangents \`a $\cal K$ sont
non pr\'ecompacts. Ainsi $\kappa$ est positive sur ${\cal K}^\prime$
 et \`a noyau de dimension $ \leq 1$. Ce noyau intersecte trivialement
 $\cal K$, car sinon, il sera un id\'eal de dimension 1 de $\cal K$,
 ce qui contredit son caract\`ere semi-simple.
\end{proof}

\section {Preuve du th\'eor\`eme alg\'ebrique}

Ce qui pr\'ec\`ede nous am\`ene \`a distinguer
 le cas  o\`u le radical $R$ est compact du  cas o\`u il
 ne l'est pas.

\subsection{Cas o\`u  le radical est compact}

Le radical \'etant compact, il est
donc ab\'elien et on a une d\'ecomposition directe:
 ${\cal G} = {\cal L}+{\cal R}$, o\`u $\cal L$ est semi-simple.
Une application comme dans la preuve pr\'ec\'edente de
la condition ($*$),
permet de montrer que $\kappa$ est d\'efinie positive sur ${\cal R}$.

Puisque  $G$ est non compact, $\cal L$
 contient  un facteur (direct) semi-simple  $\cal S$ de type non
compact.
Ainsi tout facteur de  $\cal S$ contient des vecteurs
qui d\'eterminent des
groupes \`a un param\`etre non pr\'ecompacts. Soit
${\cal S}_1$ un tel facteur. Alors,   une application
comme
dans la preuve pr\'ec\'edente  de la condition ($*$),  \`a tous les autres
facteurs de
$\cal G$ (qui centralisent
${\cal S}_1$) , permet de montrer que $\kappa$ est positive,
sur chacun de ces facteurs.
 Il s'ensuit  qu'ils sont tous compacts et
en particulier, par d\'efinition,  que $\cal S$ est simple.

Notons $\cal K$ le facteur semi-simple compact de $\cal G$.
 Le fait \ref{orthogonal} permet de montrer que la d\'ecomposition
${\cal G} = {\cal S}+{\cal K}+{\cal R}$ est orthogonale.

Montrons \`a pr\'esent que l'alg\`ebre
 simple de type non compact $\cal S$  est
 isomorphe \`a  $sl(2, {\bf R})$.

 Il est connu que dans tous les cas $\cal S$ contient une alg\`ebre
 ${\cal S}^\prime$ isomorphe \`a $sl(2, {\bf R})$. Notons
 $\cal E$ l'orthogonal  \`a ${\cal S}^\prime$. C'est
un suppl\'ementaire  de ${\cal S}^\prime$ (car ce dernier n'est pas
d\'eg\'en\'er\'e) qui est
 $ad ({\cal S}^\prime)$--invariant
 (par  bi-invariance de $\kappa$).

Il est aussi
connu (par alg\'ebricit\'e des representations d'alg\`ebres
semi-simples) que pour  $X \in {\cal S}^\prime$,
 si $ad_X$ est semi-simple (resp. nilpotent) sur
${\cal S}^\prime$, alors il en va de m\^eme
pour $ad_X$ agissant sur $\cal E$.
Il est facile de se convaincre que si tout \'el\'ement
hyperbolique (i.e.\ semi-simple \`a valeurs propres r\'eels)
$X \in {\cal S}^\prime$ agit trivialement sur
$\cal E$, alors toute l'action est triviale, et
${\cal S}^\prime$ sera un facteur direct de $\cal S$,
ce qui contredit la simplicit\'e de $\cal S$.

Par l'absurde, supposons qu'il existe $X$,  un \'el\'ement  hyperbolique
agissant
non trivialement sur
$\cal E$. Il existe donc un vecteur propre $Z \in {\cal E}$
 tel que $[X, Z] = \lambda Z$ et $\lambda \neq 0$.
Il en d\'ecoule que $Z$ d\'etermine un groupe \`a
un param\`etre non pr\'ecompact (car sinon
 $ad_Z$ serait semi-simple \`a valeurs propres
imaginaires pures). Or, il existe $Y \in {\cal S}^\prime$
nilpotent
v\'erifiant $[X, Y]= Y$. On en d\'eduit que $Z$ est
aussi vecteur propre, n\'ecessairement trivial par nilpotence
de $ad_Y$: $[Y,Z] =0$.  Donc $Y$ et $Z$
engendrent un groupe ab\'elien contenant au moins 2 groupes \`a
un param\`etre (diff\'erents) non pr\'ecompacts. Il y en
a donc un ensemble dense. Ceci contredit
l'hypoth\`ese ($*$) car $Y$ et $Z$ sont orthogonaux
et simultan\'ement isotropes. Ce dernier fait
 se voit facilement, car $\exp (ad_X)$ induit une
homoth\'etie non triviale sur chacune des directions de
$Y$ et $Z$.

Il ne reste \`a montrer du th\'eor\`eme alg\'ebrique dans notre
cas (i.e.\ lorsque $R$ est compact) que le fait que
 l'action se factorise en l'action,  d'un rev\^etement fini
de $PSL(2, {\bf R})$,  ou de mani\`ere \'equivalente
un quotient central infini de $\widetilde {SL(2, {\bf R}) }$.
 Il suffit pour cela de remarquer que
 $\widetilde {SL(2, {\bf R}) }$ ainsi que ses quotients finis,
ont tous leurs groupes \`a un param\`etres non pr\'ecompacts.
 Ce qui impliquerait que $\kappa$ est positive! \endproof

\subsection{Cas o\`u $R$ n'est pas compact}

On a alors d'apr\`es \ref{facteur.semisimple} une d\'ecomposition directe
orthogonale
 ${\cal G} ={\cal K}+{\cal R}$, o\`u ${\cal K}$ est compacte.
 Il suffit donc de montrer que $\cal R$ se d\'ecompose comme \'enonc\'e.
 On peut ainsi \`a pr\'esent oublier $\cal K$
 en supposant que $\cal G$ est r\'esoluble.

 Le nilradical
 $\cal N$ est non compact.
Consid\'erons la d\'ecomposition: ${\cal N}={\cal A} + {\cal HE}_d$
et notons $\cal Z$ le centre de
${\cal HE}_d$. Rappelons (\ref{non.canonique}) que c'est la somme
${\cal A}+{\cal Z}$ (mais pas ${\cal A}$) qui est canoniquement d\'efinie.

\subsubsection{Cas o\`u $\kappa$ n'est pas positive. Groupes
de Heisenberg tordus}

Soit $t$ un \'el\'ement de  $\cal G$ tel que $\kappa(t,t) <0$.  Il
 engendre
un groupe \`a un param\`etre    non pr\'ecompact, que l'on
 peut  supposer (apr\`es approximation) p\'eriodique,
 i.e.\ engendrant
un groupe isomorphe au cercle $S^1$.

\begin{fait}
$t$ centralise ${\cal A}+{\cal Z}$,  qui par suite engendre un
groupe (ab\'elien) compact, qui est donc en plus
 central dans $\cal G$.
\end{fait}

\begin{proof}[Preuve]  Soit $T^s = \exp(s ad_t)$ le
 groupe \`a un param\`etre d\'efini par $t$.
  Il agit sur ${\cal A}+{\cal Z}$  par transformations
orthogonales (\`a cause de la
pr\'ecompacit\'e),  en particulier semi-simples, \`a valeurs
propres de module \'egale \`a 1. Pour montrer que $t$
 centralise ${\cal A}+{\cal Z}$, il suffit de montrer que
 toute puissance non triviale $T^s$
n'a pas de sous-espace propre de dimension 2.
 Supposons par l'absurde que $\cal P$ est un tel sous-espace.
 C'est en particulier une sous-alg\`ebre de $\cal G$
car ${\cal A}+{\cal Z}$ est ab\'elienne. L'alg\`ebre
 $\cal L$ engendr\'ee
 par $t$
 et $\cal P$ est isomorphe
 \`a l'alg\`ebre de Lie du
groupe des d\'eplacements eucilidien d'un plan (engendrant
le groupe des translations-rotations du plan).

Tous les \'el\'ements de
$\cal P$ sont n\'ecessairement non pr\'ecompacts, et donc
d'apr\`es la condition ($*$), $\kappa$ est positive, non  triviale sur
$\cal P$.
Elle est donc non d\'eg\'en\'er\'ee, car son noyau est
un id\'eal propre, qui ne pourrait \^etre que $\cal P$.
En fait $\kappa$ est une forme lorentzienne
bi-invariante sur $\cal L$ (car on sait d\'ej\`a
que $\kappa(t,t) <0$).

Il suffit maintenant de remarquer qu'une telle forme, ne
 peut pas exister. En effet tout  groupe \`a un
param\`etre d\'efini par un vecteur non tangent \`a
$\cal P$ (i.e.\ qui ne soit
pas un groupe \`a un param\`etre de translations du
plan) est conjugu\'e \`a celui d\'efini par $t$,
car  c'est un groupe de rotation autour d'un certain point.
 Il s'ensuit que $\kappa$ est n\'egative en dehors de $\cal P$,
ce qui contredit son caract\`ere lorentzien.

On d\'eduit de \ref{dense.abelien} que
 ${\cal A}+{\cal Z}$ d\'etermine
un groupe compact.
\end{proof}

En fait, toujours d'apr\`es \ref{dense.abelien}, le groupe
\`a un param\`etre d\'etermin\'e par $t$
 ne commute avec aucun \'el\'ement non central de
 ${\cal HE}_d$. De plus le groupe engendr\'e par le centre
de ${\cal HE}_d$ est compact, faute de quoi, toujours d'apr\`es
 \ref{dense.abelien}, on aura $\kappa (t,t) \geq 0$.

Notons $\cal S$ l'alg\`ebre engendr\'ee par $t$  et
${\cal HE}_d$ et $S$ le groupe qu'elle d\'etermine.

Un raisonnement \'el\'ementaire
permet de voir  que $\kappa$ est lorentzienne sur
$\cal S$. On commence par constater que $\kappa$ est non
d\'eg\'en\'er\'ee, car  son noyau ne pourrait \^etre que
 le centre, et en quotientant par ce dernier, on trouve une
forme lorentzienne bi-invariante sur le produit semi-direct de $S^1$
 agissant, sans vecteur fixe, sur ${\bf C}^d$. La preuve qu'on vient
de donner ci-dessus, pour $d=1$,  de l'inexistence d'une
telle forme, se
g\'en\'eralise en toute dimension.

 Ce qui pr\'ec\'ede montre bien que
 $S$ est  un groupe de Heisenberg tordu.

Consid\'erons
 l'orthogonal ${\cal S}^\perp$. C'est bien
un suppl\'ementaire de
 $\cal S$. Par bi-invariance de $\kappa$, $[X,Y] \in
{\cal S}^\perp$ d\`es que $X \in {\cal S}$
et $Y \in {\cal S}^\perp$.  En d'autres termes, $\cal S$
 centralise le sous-espace vectoriel ${\cal S}^\perp$.
Il
en r\'esulte, puisque
${\cal HE}_d$ est un id\'eal de $\cal G$,  que
$[X,Y] =0$ d\`es que $X \in {\cal HE}_d$ et
$Y \in {\cal S}^\perp$.  Autrement dit
${\cal S}^\perp$ centralise ${\cal HE}_d$.

Soit $X \in {\cal S}^\perp$. Il centralise
${\cal N} = {\cal A}+{\cal HE}_d$, car d'apr\`es
 le fait ci-dessus $\cal A$ est central.
Il en d\'ecoule que ${\bf R} X +{\cal N}$ est une
alg\`ebre nilpotente. C'est en fait un id\'eal de
$\cal G$, car il est  connu que $[ {\cal G}, {\cal G}]
 \subset {\cal N}$ (on
avait
 suppos\'e que $\cal G$ est r\'esoluble). Par maximalit\'e de $\cal N$,
en tant qu' id\'eal nilpotent, on a: $X \in {\cal N}$.

Ainsi  ${\cal S}^\perp$ est
contenue dans le nilradical $\cal N$. On en d\'eduit
pour des raisons de dimension que ${\cal N}
={\cal S}^\perp + {\cal HE}_d$. Ainsi on peut
prendre ${\cal A}={\cal S}^\perp$. Ce qui ach\`evera la
 d\'ecomposition dans ce cas.
\endproof

\subsubsection{Cas o\`u $\kappa$ est positive}

Supposons que $\kappa$ est  positive sur $\cal G$ (suppos\'ee
r\'esoluble).
Elle admettra un noyau non trivial $\cal I$,
 sauf si   $\cal G$ est
ab\'elienne.
Supposons donc dans la suite
que $\cal I$ est non trivial.

D'apr\`es la condition ($*$), si  dim$({\cal I}) > 1$, alors
 le sous-groupe $I$ de $G$ qu'il d\'etermine est pr\'ecompact,
i.e.\ $\bar{I}$ est un tore, n\'ecessairement central. En particulier
 ${\cal I} \subset {\cal N}$. Ce qui contredit
 le fait (\ref{nilradical}) que, sur
 $\cal N$, la dimension du noyau de $\kappa$ est $\leq 1$.

Montrons que: ${\cal I} \subset {\cal N}$. En effet
sinon, ${\cal I} \cap {\cal N} = 0$. Comme $\cal I$ et $\cal N$
 sont des id\'eaux, il s'ensuit qu'ils se centraliseent l'un l'autre.
En particulier ${\cal I}+{\cal N}$ est aussi un id\'eal
nilpotent.  Ce qui contredit
  la d\'efinition de $\cal N$.
Maintenant, si $\cal G$ est nilpotente, elle se
d\'ecompose comme dans
 \ref{nilradical}. Ce qui d\'emontre le th\'eor\`eme
alg\'ebrique
 dans ce cas.
Supposons donc que $\cal G$ n'est pas nilpotente. L'alg\`ebre
quotient
est ab\'elienne car elle admet une forme
d\'efinie positive bi-invariante. Il s'ensuit que
 $[{\cal G}, {\cal G}] \subset {\cal I}$,
mais $\cal I$
n'est pas central, car sinon $\cal G$ sera nilpotente.
On en d\'eduit que si $Y$
est un g\'en\'erateur de $\cal I$,
alors le noyau de l'application $ u \to [u, Y]$ admet
un noyau $\cal L$
de codimension 1. Il existe $X$ orthogonal
\`a ${\cal L}$
 v\'erifiant $[X,Y] \neq 0$. On peut en fait supposer quitte \`a
prendre un multiple de $X$ que: $[X, Y] =Y$.
Soit ${\cal A} \subset {\cal L}$ le noyau de
$T \in {\cal L} \to [X,T] \in {\cal I}$.
Ainsi
$X$ et $Y$ engendre l'alg\`ebre de Lie du groupe affine ${\cal GA}$.
 Pour achever la preuve du th\'eor\`eme dans le pr\'esent cas,
il
suffit de montrer que $\cal A$  est une alg\`ebre centrale (elle sera
 alors imm\'ediatement un facteur direct). Comme
 par construction $\cal A$ est centralis\'e par $X$ et $Y$ et
s'injecte dans le quotient ab\'elien ${\cal G} / {\cal I}$, il suffit
juste de montrer que $\cal A$ est bien une alg\`ebre.
Soit donc $T$ et $T^\prime$ deux \'el\'ements de $\cal A$.
Par l'identit\'e de Jacobi $[X, [T, T^\prime] ]=0$. Donc
$[T, T^\prime]$ est certainement un multiple trivial de $Y$. \endproof

\section{Preuve des Th\'eor\`emes g\'eom\'etriques}

\subsection{Caract\`ere causal de l'action
lorsque $\cal G$ ne contient ni
$sl(2, {\bf R})$ ni une
alg\`ebre de Heisenberg tordue} \label{nontemporelle}

 Pour montrer que lorsque $\cal G$ ne contient pas
$sl(2, {\bf R})$ ou une alg\`ebre de
Heisenberg tordue, les orbites
sont non temporelles, il suffit d'appliquer
 \ref{fondamentale},  en remarquant  que dans ce cas, d'apr\`es
le th\'eor\`eme alg\'ebrique, les groupes \`a
 un param\`etre non pr\'ecompact sont denses.

Il s'ensuit que si $X$ est un champ isotrope au sens de $\kappa$,
alors $X(x)$ est isotrope (au
sens de $\langle \; ,\; {\rangle }_x$) pour tout $x \in M$. Pour montrer que
les orbites de $X$ sont g\'eod\'esiques, on applique le fait
suivant:
\begin{lemme} \label{orbite.geodesique} Soit $X$ un champ de Killing
 \`a norme constante: $\langle X(x),X(x)\rangle $ ne d\'epend pas de $x$.  Alors
 les orbites de $X$ sont des g\'eod\'esiques
affinement
param\'etr\'ees: $\nabla _X X (x) = 0$, pour tout $x$.
\end{lemme}
\begin{proof}[Preuve]
 En tant que champ de Killing, $X$ v\'erifie:
 $\langle \nabla_YX,X\rangle +\langle Y,\nabla_XX\rangle  =0$, pour tout champ $Y$. Mais
la constance de la	 norme entra\^{\i}ne: $\langle \nabla_YX, X\rangle  =0$.
 Par cons\'equent: $\nabla_XX=0$.
\end{proof}

\subsection{Caract\`ere localement libre des
actions des groupes de Heis\-enberg tordus}

  On supposera dans la pr\'esente section et la suivante que
 $M$ est compacte. En effet, on aura affaire dans les d\'emonstrations
 suivantes \`a des parties ferm\'ees invariantes de $M$.
La comapcit\'e de $M $ assurera l'existence de mesures invariantes
 support\'ees par ces parties. La finitude du volume de $M$ ne
l'entra\^{\i}ne  \`a priori pas. Cependant, un peu plus d'analyse de
notre situation particuli\`ere (voir \cite{Ze3}), dont on se permet de se passer
pour ne pas
encombrer davantage le texte, permet de tra\^{\i}ter ce cas l\`a.

Notre approche ressemble \`a ce niveau \`a celle
 de \cite{Gro}.

Soit $S$ un groupe de Heisenberg tordu, produit semi-direct
de $S^1$ par $He_d$,  et soit
 $Z =\{ \phi^s, s \in [0,  \pi ] \}$ son centre.
Il est facile de tirer du fait que (d'apr\`es ce qui pr\'ec\`ede)
les orbites du groupes de Heisenberg sont non temporelles, que
 les orbites de $Z$ sont isotropes.  Elles sont ainsi de plus
g\'eod\'esiques d'apr\`es le lemme \ref{orbite.geodesique}. Il s'ensuit
 que $Z$ n'admet pas de point fixe. En effet au voisinage
d'un tel point, il y aura des g\'eod\'esiques
 ferm\'ees arbitrairement petites (ce qui contredit
 la convexit\'e locale des vari\'et\'es munies de connexions affines).

Nous allons  maintenant montrer par l'absurde que l'action
 de $S$ est localement libre et ce en montrant
que sinon l'action de $Z$ ne l'est pas.
En effet soit $F$ le ferm\'e de $M$  des points
ayant un stabilisateur $S_x$ non discret. Notons ${\cal S}_x$
son alg\`ebre de Lie. On se restreint au
   ferm\'e $F_k$ o\`u la dimension de
${\cal S}_x$ est maximale \'egale \`a $k$ (certainement
$k <$ dim$(S)$ car sinon en particulier $Z$ aura
 un point fixe).
 L'action de $S$ sur $F_k$ pr\'eserve une
mesure finie $\mu$ car $F_k$ est compact
 et $S$ est r\'esoluble. La m\'ethode de preuve
suivante
est standard (voir par exemple \cite{D-G}).
 Consid\'erons
l'application de Gauss: $Ga\co F_k \to
Gr^k( {\cal S})$  qui \`a $x \in F_k$ associe
${\cal S}_x$
l'alg\`ebre de Lie
de son stabilisateur. Elle est \'equivariante
 par rapport aux actions de $S$. Ainsi
 $Ga^*( \mu)$ est une mesure sur $Gr^k({\cal S})$
 invariante par l'action de $S$.

Le lemme de Furstenberg \cite{D-G},
 s'applique aux actions des groupes alg\'ebriques. Consid\'erons
 donc la restriction de l'action pr\'ec\'edente au groupe
de Heisenberg $He_d \subset S$. D'apr\`es Furstenberg,
cette action se factorise
sur le support de la mesure, en l'action d'un
 groupe compact. Mais $He_d$ n'a aucun groupe quotient
compact non trivial.
 Il
s'ensuit  que pour $\mu$ presque tout $x$,
 ${\cal S}_x$ est normalis\'ee par
$He_d $. Si ${\cal S}_x \cap {\cal HE}_d$
est non triviale, on aura un id\'eal non trivial de
${\cal HE}_d$. Il contiendra obligatoirement le centre.
 Lorsque ${\cal S}_x$ intersecte trivialement ${\cal HE}_d$,
elle en sera
un suppl\'ementaire pour des raisons de dimension.
 Ainsi ${\cal S}_x$ sera normalis\'ee par toute l'alg\`ebre
 ${\cal S}$, c'est-\`a-dire que ${\cal S}_x$ est
un id\'eal. Ceci est impossible. \endproof

\subsection{Caract\`ere lorentzien des orbites de $S$}
\label{libre.resoluble}

Il d\'ecoule du fait que l'action
est localement libre et du fait que les orbites de
$He_d$ sont non temporelles, qu'en tout point $x$, et
pour tout $X$ tangent \`a $He_d$,
les vecteurs $X(x)$  sont de type
espace, sauf exactement celui correspondant
au centre, qui est isotrope. Pour montrer
 que  les orbites sont lorentziennes, il suffit donc
 de montrer qu'elles
sont non d\'eg\'en\'er\'ees. Or dans ce cas, le noyau
de la m\'etrique sera exactement le centre (car
l'action est localement libre).  L'ensemble des points
 \`a orbite d\'eg\'en\'er\'ee est un ferm\'e invariant. Il
supporte donc une mesure finie invariante  $\mu$. La
 forme  $L^2$ associ\'ee, i.e.\ $\kappa (X,Y) =
 \int \langle X(x),Y(x)\rangle  d \mu (x)$ est une forme
bi-invariante sur ${\cal S}$, positive
et \`a noyau exactement le centre. Ainsi le quotient
 de ${\cal S}$    par son centre admettra une m\'etrique
d\'efinie positive bi-invariante. Mais
 ceci
n'arrive pour  un groupe r\'esoluble que s'il est ab\'elien.
 \endproof

\subsection{L'orthogonal}

Notons $\cal O$ la distribution orthogonale aux
orbites
de ${\cal S}$. On va montrer que ${\cal O}
+ {\cal Z}$ est int\'egrable (o\`u $\cal Z$
est le champ de directions d\'etermin\'e
par le centre $Z$). En tout point $x$, on a une
forme antisym\'etrique:
  $ \omega: {\cal O}_x \times {\cal O}_x \to
 {\cal S}_x$, $\omega (A, B)=$ la partie
 normale du crochet $[A, B]$. L'identification
canonique de ${\cal S}_x$ \`a l'alg\`ebre de
Lie ${\cal S}$ permet d'identifier $\omega$
\`a une
forme \`a valeurs
dans ${\cal S}$. Elle v\'erifie
la relation d'\'equivariance \'evidente:
$ \omega (gA,gB) = Ad(g) \omega (A,B)$. Or la m\'etrique
 sur $\cal O$ est riemannienne, et par suite, pour
tous $A, B$ vecteurs de $\cal O$, l'orbite  $ \{ (gA,gB)
/, g \in S \}$ est pr\'ecompacte dans
 ${\cal O} \times {\cal O}$. Il en d\'ecoule
 que l'orbite de $\omega (A, B)$ par l'action adjointe
 de $S$ est pr\'ecompacte. On v\'erifie facilement
que ceci n'est le cas que du centre Donc
 $\omega $ est \`a valeurs dans $\cal Z$. ce qui
veut exactement
dire que ${\cal O} +{\cal Z}$ est int\'egrable. \endproof

Pour ce qui pr\'ec\`ede ainsi que ce qui suit, on peut consulter
respectivement \cite{Cai-Ghy} et \cite{Car-Ghy}, o\`u
l'on tra\^{\i}te de situations semblables mais
plus d\'elicates.

\subsection{Structure}
 On va transformer ``canoniquement''
la m\'etrique lorentzienne $\langle \; ,\; \rangle $
 de $M$ en une m\'etrique riemannienne
$(\;,\;)$ (qui ne sera
aucunement invariante par l'action de $S$).
 On d\'ecr\'ete que $\cal O$ reste orthogonale  aux
orbites et reste \'equip\'ee de la m\^eme m\'etrique.
 On d\'efinit la m\'etrique sur
${\cal S}_x$ par $(X(x),Y(x))=
b(X,Y)$, o\`u $b$ est un produit
scalaire d\'efini positif
quelconque (loin d'\^etre bi-invariant)
 sur ${\cal S}$.  On prendra par exemple:
$(X_i(x),X_j(x))= \delta_{ij}$ pour une certaine
base $\{ X_i \}$ de ${\cal S}$.  Le groupe
$S$ sera \'egalement \'equip\'e de la
m\'etrique
invariante \`a
 droite d\'etermin\'ee par $b$. Ainsi, pour
tout $x$, le rev\^etement $S \to  Sx$ est
isom\'etrique (cela ne veut en aucun cas dire
que $S$ agit isom\'etriquement sur l'orbite
$Sx$ au sens de la nouvelle m\'etrique
 riemannienne).

Soit $L$ la feuille du feuilletage ${\cal O}+{\cal Z}$
passant par un certain point $x_0$  munie de
la m\'etrique induite de $(\;,\;)$. Le centre $Z$ y agit
isom\'etriquement.

On a une application:
 $ p\co S  \times L \to M$,  $p(g,x) = gx$. On v\'erifie
que $p$ est une submersion riemannienne dont l'espace
horizontal est $ {\cal S} +{\cal O}$. Plus
pr\'ecis\'ement, consid\'erons $S \times _{S^1}
L$, le quotient de $S \times L$ par l'action
diagonale, et munissons le de
la m\'etrique projet\'ee de celle de ${\cal S}
 + {\cal O}$. Alors l'application induite
$\pi \co S \times_{S^1}  L \to M$ est
  localement
isom\'etrique. Par un r\'esultat bien connu
sur les applications localement isom\'etriques, $\pi$
 est un rev\^etement, car la m\'etrique
sur $S \times _{S^1} L$ est \'evidemment
compl\`ete.

Ainsi on a:  $M = \Gamma \setminus
 S \times _{S^1} L$, o\`u
$\Gamma$ est un r\'eseau de $S \times Isom_{S^1}(L)$.
Il ne reste donc du th\'eor\`eme de structure
\ref{structure.resoluble},
dans le cas des groupes de Heisenberg tordus,
qu'\`a
expliciter la m\'etrique lorentzienne sur
 $S \times _{S^1} L$. Plus pr\'ecis\'ement il s'agit de montrer que
la m\'etrique sur $S$ est essentiellement  bi-invariante
(\ref{essentiel}), ce qui fera l'objet
de la section suivante.
.

\subsection{M\'etriques lorentziennes sur $S$} \label{bi-invariante}

 On voit d'apr\`es ce qui  pr\'ec\`ede
 que la m\'etrique  lorentzienne, notons la $m$,  le long
 des orbites, qui est par hypoth\`ese invariante
par l'action \`a
gauche de $S$, doit  \'egalement \^etre
invariante \`a droite par $\Gamma_0$, la
projection
de $\Gamma$ sur $S$. Cette projection
n'est pas n\'ecessairement
discr\`ete, mais elle est \`a covolume fini, au
sens qu'il existe un sous-ensemble de
volume fini dont les it\'er\'es par
$\Gamma_0$ couvrent  $S$. Consid\'erons
 la fermeture topologique $\bar{\Gamma_0}$.
 C'est un sous-groupe unimodulaire (car
 tous les \'elements de $Ad({\cal S})$
 sont \`a valeurs propres de module
\'egale \`a 1). On peut facilement voir
que la mesure de Haar passe en une mesure
finie sur $ S / \bar{\Gamma_0}$.

 Elle d\'et\'ermine une mesure finie
 sur l'orbite de la m\'etrique $m$, invariante par l'action
 adjointe de
$S$. Donc, d'apr\`es le lemme de Furstenberg, l'action restreinte
 \`a $He_d$ se factorise en l'action d'un groupe compact. Comme ci-dessus,
ceci entra\^{\i}ne que $m$ est $Ad(He_d)$--invariante. \\

A titre de compl\'ement, on a le fait suivant qui
montre qu'il n' y a qu'une seule g\'eom\'etrie lorentzienne
 locale sur un  groupe de  Heisenberg tordu $S$. Elle m\'erite
certainement d' \^etre
mieux comprise.

\begin{proposition} Deux m\'etriques lorentziennes
bi-invariantes quelconques sur ${\cal S}$
sont \'equi\-val\-entes par un automorphisme.
\end{proposition}

\begin{proof}[Preuve]
Reprenons les notations de lapreuve de
 \ref{metrique.invariante}.  Remarquons d'abord,
qu'on peut supposer, apr\`es automorphisme, que $\beta = 0$.
 Il suffit pour cela d'appliquer un automorphisme trivial sur
 ${\cal HE}_d$ et envoyer $t$ sur $t+ \delta Z$ pour un $\delta$ convenable.

Pour normaliser le param\`etre  $\alpha$, on applique le groupe \`a
param\`etre d'homoth\'e\-ties c\'el\`ebres de l'alg\`ebre de
Heisenberg ${\cal HE}_d$. Il commute avec tous les automorphismes
et donc se prolonge trivialement  au produit semi-direct ${\cal S}$. Il
se d\'efinit ainsi: $t \to t$, $Z \to  \exp(2t) Z$ et
 $X \to \exp (t)X$,  pour $X \in {\bf C}^d$. (Ceci induit
 des homoth\'eties de $\tilde{S}$ muni de la m\'etrique
donn\'ee initialement).
\end{proof}

\subsection{Cas de $sl(2, {\bf R})$ }

D'apr\`es le th\'eor\`eme alg\'ebrique, l'action de $sl(2, {\bf R})$
s'int\'egre en
une action d'un rev\^etement fini  $PSL_k(2, { \bf R})$
de $PSL(2, {\bf R})$.  Montrons bri\`evement dans ce qui suit
le th\'eor\`eme de structure \ref{structure.semisimple} d\^u \`a \cite{Gro}.

Soit $\kappa$ la forme de Killing de $sl(2, {\bf R})$. Montrons que
si $Y \in sl(2, {\bf R})$ est isotrope au sens
de $\kappa$, alors $Y(x)$ est
isotrope au sens de $\langle \;,\; {\rangle }_x$ pour tout $x$. En effet, il est connu
 qu'un tel $Y$ est
caract\'eris\'e par le fait que $ad_X$ est nilpotent (ou de
mani\`ere \'equivalente que la matrice $2 \times 2$ et \`a trace 0,
correspondante,  est nilpotente). Il est \'egalement connu,
qu'alors il existe $X \in sl(2, {\bf R})$ tel que $[X, Y] = -Y$.
 En d'autres termes si $\phi^t$ est le flot de $X$, alors
$\phi^tY = \exp (-t) Y$. En particulier la fonction  $\langle Y,Y\rangle $
 d\'ecro\^{\i}t (exponentiellement) le long des orbites de $X$.
 Cette fonction est donc constamment nulle, car $\phi^t$ pr\'eserve
le volume.

Il en d\'ecoule que pour tout $x$, la m\'etrique restreinte \`a
l'orbite de $x$ est proportionnelle \`a $\kappa$:
$\langle X(x), Y(x)\rangle  = f(x) \kappa (X,Y)$ pour tous $X, Y$ et $x$.

Il s'ensuit en particulier qu'une orbite singuli\`ere est
isotrope. Elle est en particulier de dimension 1 ou 0.
Si elle est de dimension 1, elle sera d'apr\`es le lemme
 \ref{orbite.geodesique},
 une g\'eod\'esique (isotrope). L'action de
 $sl(2, {\bf R})$ pr\'eserve
  sa structure affine, ce qu'on  peut
facilement  voir \^etre impossible.
 L'orbite singuli\`ere est donc de dimension 0, i.e.\ un
 point fixe $x_0$
de $PSL_k(2, {\bf R})$. Soit $\phi^t$ ($t \in [0, 2 \pi]$)
 un groupe  \`a
un param\`etre de rotation de $PSL_k(2, {\bf R})$. Les
 orbites par $\phi^t$ des points proches
de $x_0$ sont des courbes
 \'egalement proches de $x_0$. On montrera dans
la suite que ces courbes
sont de type temps, i.e.\ \`a l'int\'erieur du c\^one de lumi\`ere. Ceci
est impossible, car une courbe dirig\'ee par un champ de c\^ones ne
peut pas se refermer localement.

Le point $x_0$ peut \^etre approch\'e par des points non singuliers car
l'ensemble de ces derniers points est de mesure totale ( comme pour toute
action fid\`ele pr\'eservant le volume
d'un groupe semi-simple  (voir par exemple \cite{D-G}).
 La m\'etri\-que sur l'orbite d'un point non singulier $x$ est lorentzienne,
car sinon l'orbite sera isotrope,  ce qui est
impossible car elle est de dimension
3. De plus $X(x)$ a le m\^eme caract\`ere causal que $X$, pour
tout $X \in sl(2, {\bf R})$. Comme tout $X$ engendrant
un flot non pr\'ecompact est partout non temporel, il en d\'ecoule
que les champs de type temps sont exactement ceux qui
 d\'eterminent des
flots compacts.

Enfin la m\^eme m\'ethode de preuve que pour les groupes de
Heisenberg tordus permet
de conclure que l'orthogonal est cette fois int\'egrable.
 \endproof

\subsection  {Vari\'et\'es homog\`enes}

Soit $(M, \langle \; ,\; \rangle )$ une vari\'et\'e lorentzienne  homog\`ene
 de volume fini. Son alg\`ebre de champs de Killing
 agit dessus localement transitivement. En particulier en
tout point, il y a des champs de Killing ayant un caract\`ere causal
 quelconque.
 Ceci exclut la situation d\'ecrite en \ref{nontemporelle}. En d'autres
termes, le groupe d'isom\'etries $G$ contient un
groupe $S$ qui est soit
 localement isomorphe \`a
 $SL(2, {\bf R})$, soit  isomorphe \`a
 un groupe de Heisenberg tordu. Dans
chacun des ces deux cas, d'apr\`es ce qui pr\'ec\`ede, $M$
 admet un rev\^etement qui est un produit tordu de
$S$ par une vari\'et\'e riemannienne $L$. Notre
hypoth\`ese d' homog\'en\'eit\'e nous permet de choisir $L$
 compacte. En effet comme dans les preuves pr\'ec\'edentes,
 en d\'esignant comme toujours par  $\cal O$ l'orthogonal
 aux orbites,
on prendra  pour $L$ soit une feuille de $\cal O$ lorsque
 $S$ est localement isomorphe \`a $SL(2, {\bf R})$,
 soit une feuille de ${\cal O}+{\cal Z}$ lorsque
$S$ est un groupe de Heisenberg tordu. Soit $H$ la composante
 neutre du sous-groupe de $G$ fixant (globalement) $L$.
 On d\'eduit ais\'ement du th\'eor\`eme
alg\'ebrique que $H$ est compact. Il en va de m\^eme pour
$L$, car $H$ agit transitivement dessus (\`a cause
 de l' homog\'en\'eit\'e
de $M$).

Enfin pour voir que le groupe du rev\^etement $\Gamma$ est le graphe
 d'un homomorphisme d'un r\'eseau de $S$, on remarque
 simplement que par compacit\'e de $L$, $\Gamma$
 se projette sur un groupe discret de $S$. Le noyau
de la projection de $\Gamma$ sur $S$ est un sous-groupe
fini d'isom\'etries de $L$,  qu'on peut supposer trivial en passant \`a
un quotient fini de $L$.			\endproof

\np

\renewcommand{\refname}{Bibliographie}

\Addresses\recd

\end{document}